\documentclass[a4paper,11pt]{article}
\usepackage{amsbsy,amsfonts,amsmath,amssymb,eucal,mathrsfs}
\usepackage[all]{xy}
\usepackage{pstricks,epsfig}
\usepackage{pst-plot}
\usepackage{psfrag}

\newtheorem{definition}{Definition}[section]
\newenvironment{defi}{\begin{definition} \rm}{\end{definition}}

\newtheorem{prop}[definition]{Proposition}

\newtheorem{coro}[definition]{Corollary}
\newtheorem{theo}[definition]{Theorem}
\newtheorem{notation}[definition]{Notation}

\newtheorem{construction}[definition]{Construction}

\newtheorem{remark}[definition]{Remark}
\newenvironment{rema}{\begin{remark} \rm}{\end{remark}}
\newtheorem{remarks}[definition]{Remarks}

\newtheorem{example}[definition]{Example}

\newtheorem{examples}[definition]{Examples}

\newtheorem{nothing}[definition]{$\!\!$}

\newenvironment{proo}{{\flushleft \it Proof.}}{\hfill $\square$
  \vspace{2mm}}

\newenvironment{proo-prop}{{\flushleft \it Proof of Proposition
    \ref{prop-f_4}.}}{\hfill $\square$ \vspace{2mm}}
\newenvironment{conj}{\begin{conjecture} \rm}{\end{conjecture}}
\newtheorem{conjecture}[definition]{Conjecture}

\newtheorem{definition*}{Definition}[section]
\newenvironment{defi*}{\begin{definition*} \rm}{\end{definition*}}
\newtheorem{definitions*}[definition*]{Definitions}
\newenvironment{defis*}{\begin{definitions*} \rm}{\end{definitions*}}
\newtheorem{prop*}[definition*]{Proposition}
\newtheorem{lemm*}[definition*]{Lemma}
\newtheorem{coro*}[definition*]{Corollary}
\newtheorem{theo*}[definition*]{Theorem}
\newtheorem{remark*}[definition*]{Remark}
\newenvironment{rema*}{\begin{remark*} \rm}{\end{remark*}}
\newtheorem{remarks*}[definition*]{Remarks}
\newenvironment{remas*}{\begin{remarks*} \rm}{\end{remarks*}}
\newtheorem{example*}[definition*]{Example}
\newenvironment{exam*}{\begin{example*} \rm}{\end{example*}}
\newtheorem{examples*}[definition*]{Examples}
\newenvironment{exams*}{\begin{examples*} \rm}{\end{examples*}}
\newtheorem{nothing*}[definition*]{$\!\!$}
\newenvironment{noth*}{\begin{nothing*} \rm}{\end{nothing*}}

\newtheorem{commentaire*}[definition*]{Commentaire}

\def \Xh {{{\widehat{X}}}}
\def \bp {{{\bar P}}}

\def \vs {\vskip}

\newcommand{\G}{\mathbb{G}}
\def \a {\alpha}

\newcommand{\Xt}{{\widetilde{X}}}

\newcommand{\Z}{\mathbb Z}

\def \pic {{\rm Pic}}

\newcommand{\scal}[1]{\langle #1 \rangle}

\def\cO{{\cal O}}

\begin{document}

\title{Spherical varieties and Wahl's conjecture}
\author{Nicolas Perrin}
\date{}

\maketitle

\begin{abstract}
Using the theory of spherical varieties and especially Frobenius
splitting results for symmetric varieties, we give a type independent
very short proof of Wahl's conjecture for cominuscule homogeneous
spaces for all primes different from 2.
\end{abstract}

{\def\thefootnote{\relax}
 \footnote{ \hspace{-6.8mm}
 }
 }

\centerline{\large{\bf Introduction}}

\vs 0.3 cm

Let $V$ be a smooth projective variety and let $L$ and $M$ be two line
bundles on $V$. It is natural to consider the so called \emph{Gaussian map}:
$$H^0(V\times V,{\cal I}_\Delta\otimes L\boxtimes M)\to
H^0(V,\Omega^1_V\otimes L\otimes M),$$
where $\Delta$ is the diagonal in $V\times V$, where $L\boxtimes M$ is
the external product on $V\times V$ and the map is induced by the
restriction map ${\cal I}_\Delta\to{\cal I}_\Delta/{\cal
  I}_\Delta^2\simeq\Omega^1_V$. J. Wahl studied this map in detail. In
particular in \cite{wahl} he conjectured that the Gaussian map should be
surjective for $V$ a rational homogeneous space and $L$ and $M$ any
ample line bundles. This conjecture was proved by S. Kumar in
characteristic 0 in \cite{kumar}. V. Lakshmibai, V.B. Mehta
and A.J. Parameswaran \cite{LMP} considered the situation in positive
characteristic and proved that the following conjecture (now called
\emph{LMP-conjecture}) implies Wahl's conjecture in positive
characteristic. From now on we work over an algebraically closed field $k$ of positive characteristic $p$.

\begin{conj}~

Let $V$ be a rational projective homogeneous space, let $X=V\times V$ and let $\Xt$ be the blowing-up of the diagonal $\Delta$ in $X$. Then $\Xt$ is Frobenius split compatibly with the exceptional divisor $E$.
\end{conj}

This conjecture is equivalent to the existence of a
splitting on $V\times V$ with maximal multiplicity along the
diagonal (see \cite{LMP} for more on this). This conjecture has been
considered by several authors (see for example \cite{MP}, \cite{LRS},
\cite{BL}, \cite{LT}, \cite{T}). In particular J. Brown and V. Lakshmibai in \cite{BL}
proved this conjecture for minuscule homogeneous spaces using
Representation Theoretic techniques and a case by case analysis.

In this paper we want to give a new proof of LMP-conjecture and
therefore of Wahl's conjecture for cominuscule homogeneous spaces (see
Definition \ref{def-comin}) using the theory of spherical
varieties. Let $V$ be a cominuscule homogeneous space and let $\Xt$ be
the blow-up of the diagonal in $X=V\times V$.

\begin{theo}~
\label{theo-main}

Assume that $p\neq2$, then the variety $\Xt$ is Frobenius split
compatibly with the exceptional divisor.
\end{theo}

Remark that since any minuscule homogeneous space is cominuscule
for some other group this also implies the result in the minuscule
case. The advantage of this proof is that it is mainly geometric, it
completely avoids the case by case analysis in \cite{BL} and it is very
short.

One of the main argument is to remark that if $V$ is cominuscule, then
$X=V\times V$ is spherical. Using this idea and a result of M. Brion
and S.P. Inamdar \cite{BI}, a very simple proof of Theorem
\ref{theo-main} is given in section \ref{sect-short-proof} for large
primes and in particular in characteristic 0. To obtain the result for
all odd primes, we need to do a parabolic induction from a symmetric
variety and use a result of C. de Concini and T.A. Springer
\cite{deCS}.

\paragraph{Acknowledgement} I want to thank Michel Brion for useful email
echanges on the subject in particular for the reference \cite{deCS}
and Mart\'i Lahoz Vilalta for useful discussions.

\section{Very simple proof for large primes}
\label{sect-short-proof}

In this section we give a very short proof of Theorem \ref{theo-main} for large primes using a result of P. Littelmann \cite{littelmann} and a result of M. Brion and S.P. Inamdar \cite{BI}. Write $V=G/P$ with $G$ semisimple and $P$ a parabolic subgroup and write $X=V\times V$.

\begin{theo}[Littelmann]
\label{theo-littel}~

  If $P$ is cominuscule, then the variety $X$ is spherical. In particular $\Xt$ is spherical.
\end{theo}

P. Littelmann proves that the only maximal parabolic
subgroups $P$ such that the product $G/P\times G/P$ is spherical are the minuscule
and cominuscule parabolic subgroups. In \cite{AP}, we in
particular generalise this statement and prove that this is true for
any (non-maximal) parabolic subgroup. We shall recover Theorem
\ref{theo-littel} in Section \ref{section-struct} (see Proposition \ref{prop-L/K}).

\begin{theo}[Brion-Inamdar]~

Assume that $\Xt$ is defined over $\Z$ and is spherical for some group $G$. Then  for all but finitely many primes $p$ the variety $\Xt_p$ is Frobenius split compatibly with all closed $G$-stable subvarieties.
\end{theo}

\begin{coro}~

The variety $\Xt_p$ is Frobenius split compatibly with $E_p$ for all
but finitely many primes $p$.
\end{coro}

\section{Frobenius splitting of symmetric varieties}

In this section we extend results of C. De Concini and T.A. Springer
\cite{deCS} on Frobenius splitting of compactifications of symmetric
varieties. The results we obtain are probably well know to the experts but we could not find a reference for them. We assume from now on that $p$ is not $2$ and we refer to \cite{knop} for classical results on spherical varieties.

Let $L$ be a semisimple algebraic group of adjoint type and let $K$ be
the fixed point subgroup of an involution $\theta$. Then $L/K$ is
called a homogeneous symmetric variety and it is a spherical variety
(see \cite{vust}). Let $B_L$ be a Borel subgroup of $L$ such that $B_L$ has a dense orbit in $L/K$ and let $T_L$ be a maximal torus of $L$ containing a split maximal torus $S$ (\emph{i.e.} a maximal torus such that $\theta\vert_S$ acts as the inverse). Recall the following results from \cite{deCS}.

\begin{prop}~
\label{deCS}

(\i) There exists a unique simple smooth projective toroidal compactification
${\bf Y}$ of $L/K$.

(\i\i) There exists a parabolic subgroup $Q$ of $L$ containing $B_L$, an open affine subset ${\bf Y}_0$ of ${\bf Y}$ which meets all the $L$-orbits and ${\bf
  Z}$ an affine variety contained in ${\bf Y}_0$ such that
\begin{itemize}
\item the Levi subgroup $L(Q)$ of $Q$ containing $T_L$ acts on ${\bf Z}$ and its derived subgroup $D(L(Q))$ acts trivially on ${\bf Z}$ so that ${\bf Z}$ is a toric variety for a quotient of $L(Q)/D(L(Q))$;
\item the multiplication map $R_u(Q)\times {\bf Z}\to{\bf Y_0}$ is an isomorphism.
\end{itemize}

(\i\i\i) Let $Q^-$ be the parabolic subgroup opposite to $Q$ with respect to $T_L$. Then the unique closed orbit in ${\bf Y}$ is isomorphic to $L/Q^-$ and the pull-back map $\pic({\bf Y})\to\pic(L/Q^-)$ is injective.

(\i v) The $L$-stable divisors $({\bf Y}_i)_{i\in {\bf I}}$ in ${\bf Y}$ are smooth with normal crossing and any $L$-orbit closure is the intersection of a unique
subfamily  $({\bf Y}_i)_{i\in {\bf J}}$ with ${\bf J}\subset {\bf I}$ of $L$-stable divisors.

(v) Let $\omega_{L/Q^-}$ be the canonical sheaf of the closed orbit, then it can be lifted to a line bundle ${\cal L}_{can}$ on ${\bf Y}$ and we have the formula
$$\omega_{\bf Y}={\cal L}_{can}\otimes\cO_{\bf Y}\left(-\sum_{i\in {\bf I}}{\bf Y}_i\right).$$

(v\i) Let ${\cal L}\in \pic({\bf Y})$ such that ${\cal L}\vert_{L/Q^-}$ is ample. If $k$ is even, then the restriction map $H^0({\bf Y},{\cal L}^{\otimes k})\to H^0(L/Q^-,{\cal L}^{\otimes k}\vert_{L/Q^-})$ is surjective.
\end{prop}

\begin{proo}
(\i) and (\i v) are proved in \cite[Theorem 3.9]{deCS}. (\i\i) is proved in \cite[Proposition 3.8]{deCS}. (\i\i\i) is proved in \cite[Theorem 3.9 and Theorem 4.2]{deCS}. (v) is proved in \cite[Section 5]{deCS} and (v\i) is proved in \cite[Proposition 5.7]{deCS}.
\end{proo}

In our situation, we will not deal with $L/K$ but with the quotient
$L/K^0$ where $K^0$ is the connected component of the identity in $K$
with $\vert K/K^0\vert=2$. From the former proposition we deduce the
following result on $\mathfrak{Y}$ the unique simple projective toroidal completion of $L/K^0$.

\begin{coro}~
\label{coro-deCS}

(\i) The variety $\mathfrak{Y}$ is smooth and there is a
$L$-equivariant morphism $\pi:\mathfrak{Y}\to{\bf Y}$.

(\i\i) There exists
an open affine subset $\mathfrak{Y}_0$ of $\mathfrak{Y}$ which meets all the $L$-orbits and $\mathfrak{Z}$ an affine variety contained in $\mathfrak{Y}_0$ such that
\begin{itemize}
\item the Levi subgroup $L(Q)$ of $Q$ containing $T_L$ acts on $\mathfrak{Z}$ and its derived subgroup $D(L(Q))$ acts trivially on $\mathfrak{Z}$ so that $\mathfrak{Z}$ is a toric variety for a quotient of $L(Q)/D(L(Q))$;
\item the multiplication map $R_u(Q)\times \mathfrak{Z}\to\mathfrak{Y}_0$ is an isomorphism.
\end{itemize}

(\i\i\i)
The unique closed orbit in $\mathfrak{Y}$ is isomorphic to $L/Q^-$.

(\i v) The $L$-stable divisors $(\mathfrak{Y}_i)_{i\in \mathfrak{I}}$ in $\mathfrak{Y}$ are smooth with normal crossing and any $L$-orbit closure is the intersection of a unique subfamily  $(\mathfrak{Y}_i)_{i\in \mathfrak{J}}$ with $\mathfrak{J}\subset \mathfrak{I}$ of $L$-stable divisors.

(v) The canonical sheaf $\omega_{L/Q^-}$ of the closed orbit
can be lifted to a line bundle ${\cal L}_{\mathfrak{Y}}=\pi^*{\cal L}_{can}$ on $\mathfrak{Y}$ and we have the formula
$$\omega_\mathfrak{Y}={\cal L}_{\mathfrak{Y}}\otimes\cO_\mathfrak{Y}\left(-\sum_{i\in \mathfrak{I}}\mathfrak{Y}_i\right).$$

(v\i) Let ${\cal L}\in \pic({\bf Y})$ such that ${\cal L}\vert_{L/Q^-}$ is ample. If $k$ is even, then the restriction map $H^0(\mathfrak{Y},\pi^*{\cal L}^{\otimes k})\to H^0(L/Q^-,{\cal L}^{\otimes k}\vert_{L/Q^-})$ is surjective.
\end{coro}

\begin{rema}
  The only difference with Proposition \ref{deCS} is that the restriction map on Picard
  groups is not injective any more.
\end{rema}

\begin{proo}
  (\i) The map $L/K^0\to L/K$ is the quotient by the action of the subgroup $A=K/K^0$ of the automorphism group $N_L(K^0)/K^0$ of $L/K^0$. Since $\mathfrak{Y}$ is the unique simple toroidal compactification of $L/K^0$, any automorphism of $L/K^0$ extends to an automorphism of $\mathfrak{Y}$. Taking the quotient $\mathfrak{Y}/A$, we obtain that $\mathfrak{Y}/A$ is normal. It is therefore a simple toroidal embedding of $L/K$ thus $\mathfrak{Y}/A\simeq Y$ which is smooth. In particular by \cite{zariski}, the quotient map $\pi:\mathfrak{Y}\to\mathfrak{Y}/A$ is ramified over an $L$-invariant divisor. Since these divisors are smooth and we have a double cover, the variety $\mathfrak{Y}$ is smooth.
(\i\i) Set $\mathfrak{y}=K^0/K^0$ and ${\bf y}=K/K$. Let $\mathfrak{Y}_0=\pi^{-1}({\bf Y}_0)$ which is therefore $Q$ stable, let $\mathfrak{Z}=\pi^{-1}({\bf Z})$ and let $\mathfrak{Z'}=\overline{T_L\cdot \mathfrak{y}}$ the closure of the $T_L$-orbit in $\mathfrak{Y}_0$. The map $\pi$ being $L$-equivariant we get from the isomorphism $R_u(Q)\times{\bf Z}\to {\bf Y}_0$ an isomorphism $R_u(Q)\times{\mathfrak{Z}}\to \mathfrak{Y}_0$. The equivariance of $\pi$ implies the equality $L\cdot \mathfrak{Y}_0=\mathfrak{Y}$. Furthermore, since $\pi(\mathfrak{y})={\bf y}$, the orbit $T_L\cdot \mathfrak{y}$ is contained in $\mathfrak{Z}$ which is closed so we have the inclusion $\mathfrak{Z}'\subset \mathfrak{Z}$. But the isomorphism $R_u(Q)\times\mathfrak{Z}\to \mathfrak{Y}_0$ and the fact that $\mathfrak{Y}_0$ is reduced, irreducible and normal implies that $\mathfrak{Z}$ is reduced, irreducible and normal. In particular, since $\mathfrak{Z}'$ and $\mathfrak{Z}$ have the same dimension we deduce $\mathfrak{Z}'=\mathfrak{Z}$. Therefore $\mathfrak{Z}$ is also a toric variety as claimed.
(\i\i\i) This comes from the same assertion on ${\bf Y}$.
(\i v) The $L$-orbit structure is the one of a smooth toric variety via the isomorphism given in (\i\i). The result follows.
(v) The composition $L/Q^-\to\mathfrak{Y}\to{\bf Y}$ is the inclusion of the closed orbit thus the pull-back ${\cal L}_\mathfrak{Y}=\pi^*{\cal L}_{can}$ lifts $\omega_{L/Q^-}$ to $\mathfrak{Y}$. Now on $R_u(Q)\times T_L\cdot \mathfrak{y}$, we have a nowhere vanishing volume form which is furthermore $R_u(Q)\times T_L$-invariant. In particular, there exists a canonical divisor which is a linear combination of $B_L$-stable divisors. The structure Theorem (\i\i) above and classical results on toric varieties (see for example \cite{brion-gen} and  \cite{Oda}) imply that the coefficients of the $L$-stable divisors $\mathfrak{Y}_i$ are equal to $-1$. On the other hand to compute the coefficients of the $B_L$-stable divisors which are not
$L$-stable, we only need to restrict to $L/K^0$ and since $\pi:L/K^0\to L/K$ is not ramified, these coefficients are the same as those on ${\bf Y}$. The result follows.
(v\i) The restriction map $H^0({\bf Y},{\cal L})\to H^0(L/Q^-,{\cal
  L}\vert_{L/Q^-})$ is surjective and factors through $H^0(\mathfrak{Y},\pi^*{\cal L})\to H^0(L/Q^-,{\cal
  L}\vert_{L/Q^-})$, this concludes the proof.
\end{proo}

Since ${\bf Y}$ and $\mathfrak{Y}$ are simple toroidal
and complete, for any other toroidal embedding $Y$ of $L/K$ respectively of
$L/K^0$, there exists a unique $L$-equivariant morphism of pointed
$L/K$-embeddings $\pi:(Y,y)\to({\bf Y},{\bf y})$ respectively of pointed
$L/K^0$-embeddings $\pi:(Y,y)\to(\mathfrak{Y},\mathfrak{y})$.
The
same proofs as points (\i\i) and (v) of the previous corollary give
similar results for any toroidal embedding $Y$ of $L/K$ or of $L/K^0$.
\begin{coro}~

(\i) There exists  an open subset $Y_0$ of $Y$ which meets all the $L$-orbits and $Z$ an affine toric variety contained in ${Y}_0$ such that
\begin{itemize}
\item the Levi subgroup $L(Q)$ of $Q$ containing $T_L$ acts on ${Z}$ and its derived subgroup $D(L(Q))$ acts trivially on ${Z}$ so that ${Z}$ is a toric variety for a quotient of $L(Q)/D(L(Q))$;
\item the multiplication map $R_u(Q)\times { Z}\to{Y_0}$ is an isomorphism.
\end{itemize}

(\i\i) Let $(Y_i)_{i\in I}$ be the $L$-stable divisors in $Y$ and let
$\pi:Y\to {\bf Y}$  respectively $\pi:Y\to \mathfrak{Y}$ be the
canonical $L$-equivariant morphism. We have the following formula
$$\omega_Y=\pi^*{\cal L}_{can}\otimes\cO_Y\left(-\sum_{i\in
    I}Y_i\right) \textrm{ respectively } \omega_Y=\pi^*{\cal L}_\mathfrak{Y}\otimes\cO_Y\left(-\sum_{i\in I}Y_i\right).$$
\end{coro}

Let $Y$ be any embedding of $L/K$ or of $L/K^0$, we deduce Frobenius splitting results for $Y$.

\begin{coro}~
\label{coro-split}

The variety $Y$ admits a $B_L$-canonical splitting compatible with all closed $L$-stable subvarieties.
\end{coro}

\begin{proo}
We prove this for embeddings of $L/K$, the same proof works verbatim
for embeddings of $L/K^0$. There exists a toroidal embedding $Y'$ with a birational $L$-equivariant morphism $f:Y'\to Y$ (take the normalisation of the graph of the birational transformation $Y\dasharrow{\bf Y}$ together with the projection). Furthermore any toroidal embedding has a complete toroidal embedding so we can assume (using \cite[Lemma 1.1.8]{BK}) that $Y$ is toroidal and complete. By general arguments on toric varieties and Lemma 1.1.8 in \cite{BK} again, we can also assume that $Y$ is smooth. In that case the $L$-stable divisors are also smooth and their intersection are the irreducible (smooth) closed $L$-stable subvarieties in $Y$.
Let $\tau_{L/Q^-}\in H^0(L/Q^-,\omega_{L/Q^-}^{1-p})$ be the unique
$B_L$-semi-invariant section. Since $p\neq2$ we have that $p-1$ is
even and we can use Proposition \ref{deCS} (v\i) to lift this element
to an element $\tau_{\bf Y}$ in $H^0({\bf Y},{\cal
  L}_{can}^{1-p})$. Taking the pull-back $\tau_Y=\pi^*\tau_{\bf Y}$ we
get an element in $H^0(Y,\pi^*{\cal L}_{can})$. If $(D_i)_{i\in I}$
are the $L$-stable divisors and if $(\sigma_i)$ is the canonical section of $\cO_Y(Y_i)$, multiplying with $\prod_{i\in I}\sigma_i^{p-1}$ yields an element in $H^0(Y,\omega_Y^{1-p})$. By recursive application of \cite[Exercise 1.3.E.4]{BK} and the fact that a closed $L$-orbit ${\cal O}$ in $Y$ is isomorphic to $L/Q^-$ and split by $\tau_Y\vert_{\cal O}=\tau_{L/Q^-}$ the result follows. The splitting is $B_L$-canonical since $\tau_{L/Q^-}$ and therefore $\tau_Y$ are $B_L$-semi-invariant while the sections $(\sigma_i)_{i\in I}$ are $L$-semi-invariant.
\end{proo}

\section{Structure of the open orbit}
\label{section-struct}

Let $T$ be a maximal torus of $G$ and let $W$ be the associated Weyl group. Recall that if $\varpi^\vee:\G_m\to G$ is a cocharacter of $G$ factorising through $T$, we may define a parabolic subgroup $P_{\varpi^\vee}$ of $G$ as follows:
$$P_{\varpi^\vee}=\{g\in G\ /\ \lim_{t\to0}\varpi^\vee(t)g\varpi^\vee(t)^{-1} \textrm{ exists}\}.$$
Note that $P_{\varpi^\vee}$ contains $T$. Any parabolic subgroup containing $T$ can be defined this way and there exists a unique minimal cocharacter $\varpi^\vee_P$ such that $P=P_{\varpi_P^\vee}$.

\begin{defi}~
\label{def-comin}

A parabolic subgroup is cominuscule
if its associated cocharacter $\varpi^\vee_P$ satisfies $\vert\scal{\varpi_P^\vee,\a}\vert\leq1$.
\end{defi}

Let $P$ be a parabolic subgroup, let $w_0$ be the longest element in $W$ and set $H=P\cap P^{w_0}$.

\begin{defi}~

(\i) Let $\varpi^\vee=\varpi^\vee_P+w_0(\varpi^\vee_P)$ and let $R=P_{\varpi^\vee}$ be the parabolic subgroup associated to $\varpi^\vee$.

(\i\i) Let $L$ be the Levi subgroup of
$R$ containing $T$ and let $U_R$ be the unipotent radical of $R$.

(\i\i\i) Define $(\bp,\bp')=(L\cap P,L\cap P^w)$.
\end{defi}

\begin{prop}~

(\i) The parabolic subgroup $R$ contains $H$.

(\i\i) We have the equality $\varpi^\vee_\bp+\varpi^\vee_{\bp'}=0$ as cocharacters of $L$.
\end{prop}

\begin{proo}
(\i) This is obvious by definition.
(\i\i) The roots of $L$ are the roots $\a$ such that
$\scal{\varpi^\vee_{P},\a}=-\scal{\varpi^\vee_{P^{w_0}},\a}$.
The result follows.
\end{proo}

Set $p:G/H\to G/R$ and let $K^0=\bp\cap \bp'$. We assume from now on that $P$ is cominuscule.

\begin{prop}~
\label{prop-L/K}

(\i) The morphism $p$ is locally trivial with fibre $L/K^0$.

(\i\i) There exists an involution $\theta$ of $L$ such that if
$K=L^\theta$, then $K^0$ is the connecte component of $K$ containing
the identity. Furthermore $K/K^0$ is of order 2.
\end{prop}

\begin{proo}
(\i) The map is clearly locally trivial. Its fiber is $R/H$. We need to
  check that the roots in $U_R$ are in $H$.
Such a root $\a$
  satisfies, $\scal{\varpi^\vee,\a}>0$. If $\a$ is not in $P$,
  then $\scal{\varpi_P^\vee,\a}<0$ but since
  $\scal{w_0(\varpi^\vee_P),\a}\leq1$ (since $P$ is cominuscule)
  we get a contradiction. Thus $\a$ is a root of $P$ and by the same
  argument it is a root of $P^{w_0}$ and thus of $H$.
(\i\i) The
  parabolic subgroups $\bp$ and $\bp'$ are cominuscule and opposite in
  $L$. Take $\theta$ to be the involution defined by conjugating with
  $\varpi^\vee_\bp(-1)$. Then $L^\theta=K$ and $K^0$ is the connected
  component of $K$ containing the identity. This proposition
  implies that $G/H$ and thus $X$ are spherical. Finally from the
  above description, we get easily
  get that the group $K/K^0$ is of order 2.
\end{proo}

\section{Proof of the main result}

Let $Y$ be the a $L$-equivariant embedding of $L/K^0$ and define $X_Y=G\times^RY$.

\begin{prop}~

The variety $X_Y$ is $B$-canonically split compatibly with its closed
$G$-stable subvarieties.
\end{prop}

\begin{proo}
  Choose $B'$ a Borel subgroup of $G$ such that $B'$ has a dense orbit in $X$. In particular $B'$ has a dense orbit in $X_Y$ therefore $B'_L:=B'\cap L$ has a dense orbit in $Y$. Furthermore, if $T'$ is a maximal torus of $G$ contained in $B'\cap R$ and if $R^-$ is the parabolic subgroup opposite to $R$ with respect to $T'$ we have the inclusion $B'\subset R^-$. We deduce the inclusion ${B'}^-\subset R$ where ${B'}^-$ is the Borel subgroup opposite to $B'$ with respect to $T'$.
By Corollary \ref{coro-split}, there exists a $B'_L$-canonical splitting $\varphi$ of $Y$ compatibly splitting all the $L$-stable closed subvarieties. By \cite[Proposition 4.1.10]{BK} the splitting $\varphi$ is also a ${B_L'}^-$-canonical splitting (where ${B_L'}^-={B'}^-\cap L$ is also the Borel subgroup in $L$, opposite to $B'_L$ with respect to $T'_L=T'\cap L$). Since ${B'_L}^-$ is also the quotient of ${B'}^-$ by the unipotent radical $U_R$ of $R$ the action of ${B'_L}^-$ on $Y$ induces an action of ${B'}^-$ on $Y$ and the splitting $\varphi$ is also ${B'}^-$-canonical (see for example \cite[Lemma 4.1.6]{BK}). We may therefore induce this splitting to get a ${B'}^-$-canonical splitting $\psi$ of $G\times^{{B'}^-}Y$ which splits compatibly the subvarieties $G\times^{{B'}^-}Y'$ where $Y'$ is a closed $L$-stable subvariety in $Y$ (see \cite[Proposition 4.1.17 and Exercise 4.1.E.4]{BK}). Consider the morphism $q:G\times^{{B'}^-}Y\to G\times^RY$ obtained by base change from $G/{{B'}^-}\to G/R$. We have $q_*\cO_{G\times^{{B'}^-}Y}=\cO_{G\times^{R}Y}$ therefore $\psi$ induces (see \cite[Lemma 1.1.8]{BK}) a ${B'}^-$-canonical splitting compatibly splitting the varieties $G\times^RY'$. This splitting is also a ${B'}$-canonical splitting by \cite[Proposition 4.1.10]{BK} again.
The fact that the $G$-stable closed subvarieties in $X_Y$ are of the form $G\times^RY'$ with $Y'$ any $L$-stable closed subvariety in $Y$ concludes the proof.
\end{proo}

\begin{prop}~

Any toroidal embedding of $G/H$ is of the form $G\times^RY$ where $Y$ is a toroidal embedding of $L/K^0$.
\end{prop}

\begin{proo}
  Since $G/H\simeq G\times^R L/K^0$  is induced from $L/K^0$ the lattice of $G/H$ is the same as the lattice of $L/K^0$. Now if $Y$ is a toroidal completion of $L/K^0$, the variety $G\times^R Y$ is a toroidal completion of $G/H$ (since the $G$-orbits in $G\times^R Y$ are of the form $G\times^R Y'$ with $Y'$ a $L$-orbit in $Y$). In particular the simple toroidal completion $\mathfrak{Y}$ of $L/K^0$ induces a simple toroidal completion $G\times^R\mathfrak{Y}$ of $G/H$. This in turn implies that $G/H$ and $L/K^0$ have the same valuative cone. Since the toroidal completion are classified by the non colored fans contained in the valuative cone the result follows.
\end{proo}

\begin{coro}~

The variety $\Xt$ is Frobenius split compatibly with the exceptional divisor $E$.
\end{coro}

\begin{proo}
By Proposition \ref{prop-L/K}, the variety $X$ is spherical therefore $\Xt$ is also spherical. Furthermore, any spherical variety $\Xt$ admits projective birational morphism $\Xh\to\Xt$ with $\Xh$ toroidal (take the normalisation of the graph of a birational transformation $\Xt\dasharrow X'$ where $X'$ is a projective toroidal embedding of $G/H$ -- see \cite[Lemma 5.2]{knop} for the existence of a projective toroidal embedding of $G/H$). By the former two propositions, the variety $\Xh$ is $B'$-canonically Frobenius split compatibly with its closed $G$-stable subvarieties, in particular compatibly with the inverse image of $E$. By \cite[Lemma 1.1.8]{BK} and because $\Xt$ is normal, we deduce that the variety $\Xt$ is $B'$-canonically Frobenius split compatibly with $E$.
\end{proo}

\noindent
Nicolas {\sc Perrin}, \\
{\it Hausdorff Center for Mathematics,}
Universit{\"a}t Bonn, Villa Maria, Endenicher
Allee 62,
53115 Bonn, Germany, and \\
{\it Institut de Math{\'e}matiques de Jussieu,}
Universit{\'e} Pierre et Marie Curie, Case 247, 4 place
Jussieu, 75252 Paris Cedex 05, France.

\noindent {\it email}: \texttt{nicolas.perrin@hcm.uni-bonn.de}


\begin{thebibliography}{LaRaSa090}
\bibitem[AcPe12]{AP} Achinger P., Perrin N., {\it On spherical multiple
    flags}, in preparation.
\bibitem[Bri90]{brion-gen} Brion M., {\it Vers une g{\'e}n{\'e}ralisation des espaces sym{\'e}triques}. J. Algebra {\bf 134} (1990), no. 1, 115--143.
\bibitem[BrKu05]{BK} Brion M., Kumar S., {\it Frobenius splitting methods in geometry and representation theory}. Progress in Mathematics, 231. Birkh{\"a}user, 2005.
\bibitem[BrIn94]{BI} Brion M., Inamdar S.P., {\it Frobenius splitting of spherical varieties}.  Algebraic groups and their generalizations: classical methods (University Park, PA, 1991),  207--218, Proc. Sympos. Pure Math., 56, Part 1.
\bibitem[BrLa09]{BL} Brown J. Lakshmibai V., {\it Wahl's conjecture
    for a minuscule $G/P$}, Proc. Indian Acad. Sci. Math. Sci. {\bf 119} (2009), no. 5, 571--592.
\bibitem[deCSp99]{deCS} de Concini C., Springer T.A., {\it Compactification of symmetric varieties}. Dedicated to the memory of Claude Chevalley. Transform. Groups {\bf  4}  (1999),  no. 2-3, 273--300.
\bibitem[Kno89]{knop} Knop F., {\it The Luna-Vust theory of spherical embeddings}. Proceedings of the Hyderabad Conference on Algebraic Groups (Hyderabad, 1989), 225--249.
\bibitem[Kum92]{kumar} Kumar S., {\it Proof of Wahl's conjecture on
    surjectivity of the Gaussian map for flag va\-rie\-ties},
  Amer. J. Math. {\bf 114} (1992) , no. 6, 1201--1220.
\bibitem[LaMePa98]{LMP} Lakshmibai V., Mehta V.B., Parameswaran A.J.,
  {\it Frobenius splitting and blow-ups}, J. Algebra {\bf 208} (1998),
  no. 1, 101--128.
\bibitem[LaRaSa09]{LRS} Lakshmibai V., Raghavan K.N., Sankaran P.,
  {\it Wahl's conjecture holds in odd characteristic for symplectic
    and orthogonal Grassmannians}, Cent. Eur. J. Math. {\bf 7} (2009),
  no. 2, 214--223.
\bibitem[LaTh10]{LT} Lauritzen N., Thomsen J.F., {\it Maximal compatible splitting and diagonals of Kempf varieties }, to appear in Annales de l'Institut Fourier, arXiv:1004.2847.
\bibitem[Lit94]{littelmann} Littelmann P., {\it On spherical double cones}.  J. Algebra {\bf 166} (1994),  no. 1, 142--157.
\bibitem[MePa97]{MP} Mehta V.B., Parameswaran A.J., {\it On Wahl's
    conjecture for the Grassmannians in positive characteristic},
  Intern. J. Math. {\bf 8} (1997), no. 4, 495--498.
\bibitem[Oda88]{Oda} Oda T., {\it Convex bodies and algebraic geometry}.
An introduction to the theory of toric varieties. Ergebnisse der Mathematik und ihrer Grenzgebiete (3), 15. Springer-Verlag, Berlin, 1988.
\bibitem[Tho10]{T} Thomsen J.F., {\it A proof of Wahl's conjecture in the symplectic case}, arXiv:1009.0479.
\bibitem[Vus74]{vust} Vust T., {\it Op{\'e}ration de groupes r{\'e}ductifs dans un type de c{\^o}nes presque homog{\`e}nes}. Bull. Soc. Math. France {\bf 102} (1974), 317--333.
\bibitem[Wah91]{wahl} Wahl J., {\it Gaussian maps and tensor products of
    irreducible representations}, Manu\-scripta Math. {\bf 73} (1991),
  no. 3, 229--259.
\bibitem[Zar58]{zariski} Zariski O., {\it On the purity of the branch locus of algebraic functions.} Proc. Nat. Acad. Sci. U.S.A. {\bf 44} (1958) 791{\^a}€"796.
\end{thebibliography}
\end{document}